%% file: toric-fano.tex
\documentclass[oneside]{amsart}
\usepackage{etex}
\usepackage{pstricks}
\usepackage{amsmath}
\usepackage{amsfonts}
\usepackage{amsthm}
\usepackage{hyperref}
\usepackage{amscd}
\usepackage{color}
\usepackage{blkarray}
\usepackage{enumerate}
\usepackage{subcaption}
\usepackage{xypic}
\usepackage{tikz}
\DeclareCaptionSubType[Alph]{figure}

\newcommand{\bF}{\mathbf{F}}

\newcommand{\ZZ}{\mathbb{Z}}
\newcommand{\KK}{\mathbb{K}}
\newcommand{\PP}{\mathbb{P}}

\newcommand{\U}{\mathcal{U}}
\newcommand{\A}{\mathcal{A}}

\newcommand{\Aff}{\mathbb{A}}
\newcommand{\RR}{\mathbb{R}}
\newcommand{\QQ}{\mathbb{Q}}
\newcommand{\NN}{\mathbb{N}}

\newcommand{\co}{\mathfrak{o}}
\newcommand{\mcL}{\mathcal{L}}

\newcommand{\Gr}{\mathbb{G}}
\newcommand{\mgeq}{\succeq}

\DeclareMathOperator{\Bl}{Bl}

\DeclareMathOperator{\spec}{Spec}
\DeclareMathOperator{\proj}{Proj}

\newtheorem*{mainthm}{Main Theorem}
\newtheorem{thm}{Theorem}[section]
\newtheorem{prop}[thm]{Proposition}
\newtheorem{lemma}[thm]{Lemma}
\newtheorem{cor}[thm]{Corollary}

\theoremstyle{definition}

\newtheorem{defn}[thm]{Definition}
\newtheorem{rem}[thm]{Remark}

\newtheorem{ex}[thm]{Example}

\title{On Fano Schemes of Toric Varieties}

\author{Nathan Ilten}
\address{Department of Mathematics, Simon Fraser University,
8888 University Drive, Burnaby BC V5A1S6, Canada}
\email{nilten@sfu.ca}

\author{Alexandre Zotine}
\address{Department of Mathematics, Simon Fraser University,
8888 University Drive, Burnaby BC V5A1S6, Canada}
\email{szotine@sfu.ca}

\include{figures}
\begin{document}

\begin{abstract}
Let $X_\A$ be the projective toric variety corresponding to a finite set of lattice points $\A$. We show that irreducible components of the Fano scheme $\bF_k(X_\A)$ parametrizing $k$-dimensional linear subspaces of $X_\A$ are in bijection to so-called \emph{maximal Cayley structures} for $\A$. We explicitly describe these irreducible components and their intersection behaviour, characterize when $\bF_k(X_\A)$ is connected, and prove that if $X_\A$ is smooth in dimension $k$, then every component of $\bF_k(X_\A)$ is smooth in its reduced structure. Furthermore, in the special case $k=\dim X_\A-1$, we describe the non-reduced structure of $\bF_k(X_\A)$. Our main result is closely related to concurrent work done independently by Furukawa and Ito.
\end{abstract}
\maketitle
\section{Introduction}
Given an embedded projective variety $X\subset \PP^n$, its $k$th Fano scheme $\bF_k(X)$ is the fine moduli space parametrizing $k$-dimensional projective linear subspaces of $X$. Such Fano schemes have been extensively studied when $X$ is a generic (or at least smooth) hypersurface, see e.g.  \cite{altman:77a,barth:81a,langer:97a,harris:98a}, or more generally when $X$ is a complete intersection \cite{debarre:98a}. Only several isolated cases of  Fano schemes have been studied for varieties $X$ which are not generic (or at least smooth) complete intersections, see for example \cite{ilten:15a,larsen:14a}. 

In this article, we consider the Fano scheme $\bF_k(X_\A)$ when $X_\A$ is the projective  toric variety associated to a finite set of lattice points  $\A$, see \S\ref{sec:toric}. Our main result is a  complete description of the irreducible components of $\bF_k(X_\A)$ in their reduced structure, which we now summarize. A \emph{face} of $\A$ is the intersection of $\A$ with a face of its convex hull.
A \emph{Cayley structure} for such a face $\tau$ is an affine-linear surjective map $\pi$ from $\tau$ to the set $\Delta_l$ of standard basis vectors in $\ZZ^{l+1}$. In \S \ref{sec:cayley} we define a natural partial order on the set of Cayley structures for faces of $\A$. 
\begin{mainthm}[Theorem \ref{thm:main}]
There is a bijection between irreducible components of $\bF_k(X_\A)$ and maximal Cayley structures $\pi:\tau\to\Delta_l$ where $\tau$ is a face of $\A$ and $l\geq k$.
\end{mainthm} 
We give an explicit description of the component $Z_{\pi,k}$ of $\bF_k(X_\A)$ corresponding to a Cayley structure $\pi$ in \S\ref{sec:comps}. These components are locally toric, and if a general $k$-plane parametrized by $Z_{\pi,k}$ does not extend to a $k+1$-plane in $X_\A$, the component is (globally) a toric variety. We say that such components are \emph{components of maximal $k$-planes}; they correspond to maximal Cayley structures $\pi:\tau\to \Delta_l$ for $l=k$.
In Theorem \ref{thm:polytope}, we describe the set of lattice points $\A_{\pi}$ such that a component $Z_{\pi,k}$ of maximal $k$-planes is given by the projective toric variety $X_{\A_\pi}$ with respect to the Pl\"ucker embedding of $\bF_k(X_A)$. 
We also give a simple formula for the dimension of all components $Z_{\pi,k}$ (Proposition \ref{prop:dim}).

Our explicit understanding of the irreducible components of $\bF_k(X_\A)$ allows us to give a combinatorial characterization of the connected components of $\bF_k(X_\A)$ (Theorem \ref{thm:connected}). Furthermore, we show that if $X_\A$ is non-singular in dimension $k$, then every component of $\bF_k(X_\A)$ is smooth (in its reduced structure), see Corollary \ref{cor:regular}. 
Finally, we are able to explicitly describe the intersection behaviour of the components (Theorem \ref{thm:intersect}).

In \cite{casagrande:08a}, Casagrande and Di Rocco show that a $\QQ$-factorial normal projectively embedded toric variety $X_\A$ is covered by lines if and only if there is a Cayley structure $\pi:\A\to \Delta_1$. Ito generalizes this in \cite{ito:15a} to show that a polarized normal toric variety $(X,\mcL)$ is covered by $k$-planes if and only if the set of lattice points $\A$ corresponding to eigensections of $\mcL$ admits a Cayley structure $\pi:\A\to \Delta_k$. As a corollary of our main theorem, we recover Ito's result in the case that $\mcL$ defines an embedding (i.e. is very ample). While Ito uses an ingenious degeneration argument to produce a Cayley structure, our argument constructs the Cayley structure directly by considering the columns of a $(k+1)\times\#\A$ matrix whose rowspan is a $k$-dimensional linear space of $X_\A$ intersecting the torus, see Remark \ref{rem:cover}.

The key step in proving our main theorem (Theorem \ref{thm:main}) is to prove that for any linear subspace $L\subset X_\A$, there is a corresponding Cayley structure $\pi$ such that $[L]$ is contained in $Z_{\pi,k}$. This statement has been proven concurrently by Furukawa and Ito in independent work using a different argument, see \cite[Theorem 3.2]{furukawa:16a}. They make use of this result to describe the dual defect of a toric variety. 

As a first step towards understanding the (potentially) non-reduced structure of $\bF_k(X_\A)$, we locally give a combinatorial description of the scheme structure of  $\bF_{k}(X_\A)$  when $k=\dim X_\A-1$ and $X_\A$ is smooth in codimension one, see \S\ref{sec:mult}. In particular, we give formulas for the multiplicities of isolated hyperplanes in $X_\A$. This generalizes the results of the first author concerning the Fano scheme of lines on a toric surface \cite{ilten:14a}. This codimension-one case is particularly amenable to study since in this situation, each toric fixed point of $\bF_{k}(X_\A)$ is contained in a single irreducible component.

Our study of Fano schemes is partially motivated by a number of recent applications. Kiraly and Larsen use Fano schemes to solve an identifiability problem in stationary subspace analysis \cite{larsen:14a}. Results on the Fano schemes for determinantal and permanental hypersurfaces have found applications in geometric complexity theory, leading to new components of the boundary of the orbit closure of the determinant \cite[\S 5]{landsberg} and to a quadratic lower bound on the determinantal complexity of the permanent \cite{landsberg:13a}. Furthermore, recent work by the first author and Teitler used results on certain Fano schemes to determine the product or Chow rank of the $3\times 3$ permanent and determinant polynomials \cite{ilten:16a}. In fact, the key calculation of \cite[Proposition 3.2]{ilten:16a} is closely related to the toric case, see Remark \ref{rem:chow}. We hope that our results here may lead to further similar applications.

\subsection*{Acknowledgements} The first author was partially supported by an NSERC discovery grant. The second author was partially supported by an SFU-VPR undergraduate student research award. We thank Dustin Cartwright, Christian Haase, and Atsushi Ito for helpful conversations.

\section{Background and Notation}\label{sec:background}
\subsection{Toric Varieties}\label{sec:toric}
For details on toric varieties, see e.g.~\cite{CLS}. 
Throughout this article, we work over an algebraically closed field $\KK$.
Let $M$ be lattice, and $\A$ a finite set of elements of $M$. 
By $M_\A$ we denote the lattice generated by $u-v$ for $u,v\in \A$, and by $S_\A$ we denote the semigroup in $M\times\ZZ$ generated by $(u,1)$ for $u\in\A$.
 The projective toric variety associated to $\A$ is 
\[
X_\A=\proj \KK[S_\A ].
\]
Here, $\KK[S_\A]$ is the semigroup algebra associated to $S_\A$, and the $\ZZ$ grading of this algebra is given by projection onto the final $\ZZ$-factor.
More concretely, $X_\A$ is the subvariety of the projective space $\PP^{\#\A-1}$ with coordinates $x_u$ for $u\in \A$ cut out by binomial equations of the form
\[
\prod x_u^{a_u}=\prod x_u^{b_u}
\] 
for $a_u,b_u\in\ZZ_{\geq 0}$ satisfying the affine relation
\[
\sum a_uu=\sum b_uu,\qquad \sum a_u=\sum b_u.
\]

We denote by $\chi^w$ the element of $\KK[S_\A]$ corresponding to $w\in S_\A$.
The \emph{dimension} of $\A$ is the dimension of the smallest affine subspace of $M\otimes \RR$ containing it, which is equivalently the rank of $M_\A$ and the dimension of $X_\A$.
The variety $X_\A$ comes equipped with a faithful action by the torus $T=\spec \KK[M_\A]$ induced by the $M\times \ZZ$-grading of $\KK[S_\A]$. 

A \emph{face} $\tau$ of $\A$ is a set of the form $\A\cap F$, where $F$ is a face of the convex hull of $\A$. If $\tau$ is a face of $\A$, we write $\tau\prec \A$. For $\tau\prec \A$, there is a natural closed embedding
\[
X_\tau\hookrightarrow X_\A
\]
given by setting $\chi^{(u,1)}=0$ for $u\notin \tau$.
This induces a bijection between faces of $\A$ and orbits of $X_\A$: the orbit corresponding to a face $\tau\prec \A$ is given by the orbit of a general point of $X_\tau\subset X_\A$ \cite[Corollary 3.A.6]{CLS}.

\begin{figure}
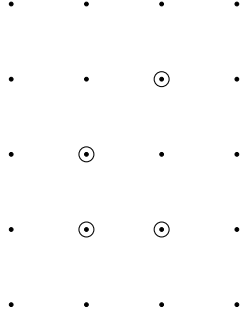

\lppic
\caption{Lattice points in Example \ref{ex:surface1}\label{fig:ex}}
\end{figure}

\begin{ex}\label{ex:surface1}
Consider the set $\A$ whose elements in $\ZZ^2$ are given by the columns of
\[
\left(\begin{array}{c c c c}
0& 0& 1& 1\\
0& 1& 0& 2
\end{array}\right).
\]
This collection of lattice points is pictured in Figure \ref{fig:ex}. The corresponding projective variety is a non-normal hypersurface cut out by the single binomial $xy^2=zw^2$. $\A$ has four $1$-dimensional faces, all of which are empty simplices.
\end{ex}

\begin{ex}\label{ex:birkhoff1}
Let $\A$ be the subset of $\ZZ^{3\times 3}$ consisting of the six $3\times 3$ permutation matrices. Then $\A$ is $4$-dimensional, and its convex hull is known as the Birkhoff polytope $B_3$. If we label the elements of $\A$ by $u_0,u_1,u_2,v_0,v_1,v_2$ with $u_i$ representing an even permutation and $v_j$ representing an odd one, the projective variety $X_{\A}$ is cut out by the single binomial
\[
x_{u_0}x_{u_1}x_{u_2}=x_{v_0}x_{v_1}x_{v_2}.
\]
The set $\A$ has exactly nine $3$-dimensional faces, given by omitting exactly one of the $u_i$ and one of the $v_j$.
\end{ex}

\subsection{Fano Schemes}\label{sec:fano}
Let  $X\subset \PP^n$ any projective $\KK$-scheme. For any natural number $k\in \NN$, the \emph{$k$th Fano scheme} $\bF_k(X)$ is the fine moduli space parametrizing $k$-planes of $\PP^n$ contained in $X$, see e.g. \cite[\S IV.3]{eisenbud:00a} for a detailed description. The scheme $\bF_k(X)$ is a subscheme of the Grassmannian $\Gr(k,n)$ parametrizing $k$-planes of $\PP^n$. For all of this article except for \S\ref{sec:mult}, the reader may forget the scheme structure of $\bF_k(X)$ and simply think of it as the subvariety of $\Gr(k,n)$ whose points correspond to $k$-planes of $\PP^n$ contained in $X$.
Given a $k$-plane $L\subset X$, we denote the corresponding point of $\bF_k(X)$ by $[L]$.

The Pl\"ucker embedding of $\Gr(k,n)$ induces an embedding of $\bF_k(X)$ in $\PP^{{n+1}\choose {k+1}}$. We will call the corresponding affine charts of $\bF_k(X)$ the \emph{Pl\"ucker charts}.
Concretely, these Pl\"ucker charts may be thought of as follows. Any $k$-plane of $\PP^n$ may be represented non-uniquely as the rowspan of a $(k+1)\times(n+1)$ matrix $P$ of full rank. The ${n+1}\choose{k+1}$ Pl\"ucker charts are obtained by choosing some $(k+1)\times(k+1)$ square submatrix and imposing the condition that it be invertible. On this chart, any $k$-plane may be  represented uniquely as the rowspan of a $(k+1)\times(n+1)$ matrix, where we impose the condition that the chosen square submatrix is just the identity matrix. Coordinate functions on this chart are given by the remaining entries of the matrix.

 If $G$ is an algebraic group acting on $\PP^n$ which fixes $X$, then the Fano schemes $\bF_k(X)$ inherit a  $G$-action. In particular, the Fano scheme $\bF_k(\A)$ has a natural torus action. 
The torus fixed points of the Fano schemes $\bF_k(X_\A)$ are easy to describe. By an \emph{empty $k$-simplex} we mean $k+1$ lattice points whose affine span is $k$-dimensional.
\begin{prop}\label{prop:fixed}
The torus fixed points of $\bF_k(X_\A)$ are in bijection with the set of faces of $\A$ which are empty $k$-simplices. 
\end{prop}
\begin{proof}
Any fixed point of $\bF_k(X_\A)$ must be the closure of a $k$-dimensional $T$-orbit in $X_\A$, which are in bijection to the set of $k$-dimensional faces of $\A$. Given a face $\tau\subset \A$, the corresponding orbit closure as a subvariety of $\PP^n$ is simply $X_\tau$ sitting inside of a linear subspace of $\PP^n$. But $X_\tau$ is a linear space if and only if there are no affine relations among the elements of $\tau$, that is, $\tau$ is an empty $k$-simplex.
\end{proof}
\noindent We denote by $L_\sigma$ the $k$-plane of $X_\A$ corresponding to an empty $k$-simplex $\sigma$.

\section{Cayley Structures and Main Result}\label{sec:cayley}
We denote by $\Delta_l\subset \ZZ^{l+1}$ the standard basis vectors $e_0,\ldots,e_{l}$.
 Fix a finite set of lattice points $\A\subset M$.
\begin{defn}
Let $\tau$ be a face of $\A$.
A \emph{Cayley structure} for $\tau$ is a surjective map $\pi:\tau\to \Delta_l$ which preserves affine relations.
\end{defn}

\begin{rem}
If $\A$ is the set of lattice points of a lattice polytope $P$, then the existence of a Cayley structure $\pi:\A\to\Delta_l$ means that $P$ can be written as a Cayley polytope of length $l+1$.
\end{rem}

\begin{figure}
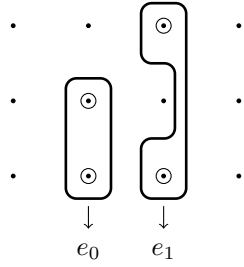

\cspic
\caption{A Cayley Structure in Example \ref{ex:surface2}\label{fig:ex2}}
\end{figure}

\begin{ex}\label{ex:surface2}
Considering the set $\A$ from Example \ref{ex:surface1}, the map $\A\to \Delta_1$ sending $\{(0,0),(0,1)\}$ to $e_0$ and $\{(1,0),(1,2)\}$ to $e_1$ is a Cayley structure, see Figure \ref{fig:ex2}. On the other hand, the map $\A\to \Delta_1$ sending $\{(0,0),(1,0)\}$ to $e_0$ and $\{(0,1),(1,2)\}$ to $e_1$ is not a Cayley structure, since it does not preserve the affine relation \[2\cdot (0,0)+(1,2)=(1,0)+2\cdot(0,1).\]
\end{ex}

A Cayley structure on a face $\tau$ of $\A$ determines a toric subvariety $Z_\pi$ of $\bF_l(X_\A)$ as follows. Indeed, since $\tau$ is a face of $\A$, we have a closed embedding $X_\tau\hookrightarrow X_\A$, see \S\ref{sec:toric}. On the other hand, a Cayley structure $\pi:\tau\to \Delta_l$ induces a surjection
\begin{align*}
\KK[S_\tau]&\to \KK[y_{e_0},\ldots,y_{e_l}]\\
\chi^{(u,1)}&\mapsto y_{\pi(u)}
\end{align*}
leading to a closed embedding $\PP^l\hookrightarrow X_\tau$, and thus an $l$-plane $L_\pi\subset X_\A$ and corresponding point $[L_{\pi}]\in \bF_l(X_\A)$. The variety $Z_\pi$ is defined to be the torus orbit closure of $[L_{\pi}]$ in $\bF_k(X_\A)$. By construction, it is a (potentially non-normal) toric variety.

Given a Cayley structure $\pi:\tau\to \Delta_l$, we also construct a subvariety $Z_{\pi,k}$ of $\bF_k(X_\A)$ for any $k\leq l$. Indeed, we take $Z_{\pi,k}$ to be the subvariety of $\bF_k(X_\A)$ whose points correspond to those $k$-dimensional linear spaces which are contained in any $l$-dimensional linear space $L'$ parameterized by $Z_\pi$. This set $Z_{\pi,k}\subset \bF_k(X_\A)$ is indeed a variety: consider the pullback $\U$ to $Z_\pi$ of the universal bundle on $\bF_l(X_\A)$. Then the Grassmann bundle $\Gr(k,\U)$ parametrizing projective $k$-planes in the fibers of $\U$ has a natural proper map to $\bF_k(X_\A)$. Its image is closed, and it exactly $Z_{\pi,k}$ as described above.

In order to state our main theorem, we define a partial order on the set of Cayley structures. We define $\pi:\tau \to \Delta_l$ to be greater than or equal to $\pi':\tau'\to \Delta_{l'}$ (written $\pi\mgeq \pi'$) if 
\begin{enumerate}
\item $\tau'$ is a face of $\tau$; and 
\item There is a surjection $\rho:\Delta_{l}\to\Delta_{l'}$ such that $(\rho\circ\pi)_{|\tau'}=\pi'$. 
\end{enumerate}
In other words, $\pi\mgeq\pi'$ if there exists a map $\rho$ making the diagram commute:
\[
\xymatrix{
\tau' \ar[r]^{\pi'}\ar@{^{(}->}[d]  &\Delta_{l'}\\
\tau \ar[r]^{\pi} & \Delta_{l} \ar@{-->}[u]^{\rho}}.
\]
Note that when considering Cayley structures $\pi:\tau\to\Delta_l$ and $\pi':\tau\to\Delta_l$, we do not differentiate between $\pi$ and $\pi'$ if $\pi\mgeq \pi'\mgeq \pi$, that is, they differ by a permutation of $\Delta_l$. We call such Cayley structures equivalent.

\begin{thm}\label{thm:main}
For $k\geq 1$, the irreducible components of $\bF_k(X_\A)$ with their reduced structure are exactly
the varieties $Z_{\pi,k}$, as $\pi$ ranges over all maximal Cayley structures $\pi:\tau\to \Delta_l$ for $\tau\prec \A$ and $l\geq k$. 
\end{thm}
\noindent We will prove this theorem in \S\ref{sec:proof}. First, we need to better understand the varieties $Z_{\pi,k}$, which we do in the next section.

\begin{ex}\label{ex:birkhoff2}
Consider the set $\A$ from Example \ref{ex:birkhoff1} whose convex hull is the Birkhoff polytope $B_3$. There are exactly six Cayley structures (up to equivalence) of the form $\pi:\A\to\Delta_2$, given by by projecting a permutation matrix in $\ZZ^{3\times 3}$ onto a single row or a single column. These Cayley structures are all maximal. 
On the other hand, the nine facets $\tau\prec \A$ described in Example \ref{ex:birkhoff1} are all empty $3$-simplices, so they naturally give Cayley structures $\pi:\tau\to\Delta_2$. These Cayley structures are also maximal. In fact, every maximal Cayley structure for $\A$ is of the above form. By Theorem \ref{thm:main}, we conclude that $\bF_k(X_\A)$ has $15$ irreducible components for $k=1,2$ and $\bF_3(X_\A)$ has $9$ irreducible components.
\end{ex}

\section{Geometry of $Z_{\pi,k}$}\label{sec:comps}
Let $\pi:\tau\to \Delta_l$ be a Cayley structure.
The $l$-plane $L_\pi$ arises as the rowspan of the $(l+1)\times \# \A$ matrix $P=(p_{iu})$ where 
\[
p_{iu}=
\begin{cases}
1 & u\in \tau \ \textrm{and}\ \pi(u)=e_i\\
0 & \textrm{otherwise}.
\end{cases}
\]
Here, the rows of $P$ are indexed by $i=0,\ldots,l$ and the columns by $u\in \A$.
Given $t\in T$, we then have that $t\cdot L_\pi$ is the rowspan of the matrix $t\cdot P$ where $T$ acts on column $u$ by multiplication with $\chi^u$, that is, 
\[
t\cdot P=(\chi^u(t) p_{iu}).
\] 

We now describe the local structure of $Z_\pi$ on the Pl\"ucker charts of $\bF_l(X_\A)$.
Choose some subset $\sigma\subset \A$ of size $l+1$ for which the corresponding minor of $P$ is non-vanishing. This occurs if and only if $\sigma\subset \tau$ and $\sigma$ contains exactly one element from each fiber $\pi^{-1}(e_i)$. 
\begin{defn}Fixing $\sigma\subset\tau$ as above, for each $u\in \tau$ we define 
$\lambda(u)$ to be the unique $v\in\sigma$ such that $\pi(v)=\pi(u)$.
\end{defn}
Setting the Pl\"ucker coordinate corresponding to $\sigma$ equal to one,
$t\cdot L_\pi$ is expressed uniquely as the rowspan of the matrix $(p_{vu}(t))$, with 
\begin{align*}
p_{vu}(t)=\begin{cases}
\chi^{u-\lambda(u)}(t)& u\in \tau\ \textrm{and}\ \pi(u)=\pi(v)\\
0 & \textrm{else}.
\end{cases}
\end{align*}
Here, the rows and columns of $(p_{vu}(t))$ are labeled respectively by elements of $\sigma$ and $\A$.

On the Pl\"ucker chart corresponding to $\sigma$,  $Z_\pi$ is thus isomorphic to the affine toric variety whose coordinate ring is
\[
\KK[S(\pi,\sigma)]
\]
where $S(\pi,\sigma)$ is the semigroup generated by the lattice elements $u-\lambda(u)$ for $u\in\tau$.
If $\sigma$ is a face of $\A$, then the semigroup $S(\pi,\sigma)$ is pointed. Indeed, let $w\in M^*$ be such that $\sigma$ is the face of $\A$ on which $w$ is minimized. Then $w(u-\lambda(u))\geq 0$, with equality if and only if $u\in\sigma$. Hence, every non-trivial element $s$ of $S(\pi,\sigma)$ satisfies $w(s)>0$, so $S(\pi,\sigma)$ is pointed.  
In such situations, this chart of $Z_\pi$ has a unique torus fixed point \cite[Proposition 1.3.2]{CLS}, which is exactly the fixed point corresponding to $\sigma$ under Proposition \ref{prop:fixed}. We thus have
\begin{prop}
The torus fixed points of $Z_\pi$ are in bijection with empty simplicial $l$-faces of $\tau$ surjecting onto the vertices of $\Delta_l$. 
\end{prop}

We also need to understand the local structure of $Z_{\pi,k}$ for a Cayley structure $\pi:\tau\to \Delta_l$, $l\geq k$.
Similar to above, the Pl\"ucker charts containing $Z_{\pi,k}$ correspond to choosing some subset $\sigma\subset \A$ of size $k+1$ for which the corresponding submatrix of $P$ has full rank. This  is the same as requiring that $\sigma$ surjects onto $k+1$ elements of $\Delta_l$.

To locally describe $Z_{\pi,k}$, we first restrict to a local chart of $Z_\pi$ corresponding to some $\widetilde\sigma\subset \tau$ as above with $\widetilde\sigma$ in bijection with the vertices of $\Delta_l$. On this chart of $Z_\pi$, the universal bundle $\U$ trivializes by means of the matrix $(p_{vu}(t))$ from above. The projective $k$-planes parametrized by the restriction of $\Gr(k,\U)$ to this chart are thus the $k$-dimensional subspaces of the $l$-dimensional rowspan of $(p_{vu}(t))$ in projective space. The local Pl\"ucker charts of this relative Grassmannian correspond to choosing a $k$-face $\sigma$ of $\widetilde\sigma$. On such a chart, a parametrization of the corresponding linear spaces is given by the rowspan of the matrix $(q_{vu}(t))$ where 
\[
q_{vu}(t)=\begin{cases}
\chi^{u-\lambda(u)}(t) & u\in \tau\ \textrm{and}\ v=\lambda(u)\in\sigma\\
\lambda_{v\lambda(u)}\chi^{u-\lambda(u)}(t) & u\in \tau\ \textrm{and}\ \lambda(u)\notin\sigma\\
0 & \textrm{otherwise}
\end{cases}
\]
for parameters $\lambda_{vw}$ with $v\in\sigma$, $w\in \widetilde\sigma\setminus\sigma$.
This matrix looks like
\begin{equation}\label{eqn:matrix}
\begin{blockarray}{cccccc}
& \lambda(u)=v_0 &  \lambda(u)=v_1 &  \lambda(u)=v_k & \lambda(u)=w\notin\sigma' & u\notin \tau \\
\begin{block}{c(c|c|c|c|c)}
v_0   &  \chi^{u-v_0}  &0  &0&\lambda_{v_0w}\chi^{u-w} & 0  \\
v_1   & 0&  \chi^{u-v_{1}}    &0&\lambda_{v_1w}\chi^{u-w} & 0  \\
  & \vdots       & &&\vdots & \vdots  \\
v_k &  0&0&  \chi^{u-v_{k}} &\lambda_{v_kw}\chi^{u-w} & 0  \\
\end{block}
\end{blockarray}
\end{equation}
where $\sigma=\{v_{0},\ldots,v_{k}\}$. 
Note that to obtain this matrix, we have simply taken the $k+1$ rows of $(p_{vu}(t))$  corresponding to the elements of $\sigma$, and then added on generic linear combinations of the remaining rows.

Let $S(\pi,\widetilde\sigma,\sigma)$ be the semigroup of $M\times \ZZ^{\sigma\times(\widetilde\sigma\setminus\sigma)}$ generated by the lattice elements
\begin{equation}\label{eqn:gens}
\begin{aligned}
\gamma(u):=&u-\lambda(u) &\qquad u\in \tau\ \textrm{and}\ \lambda(u)\in\sigma\\
\gamma(v,u):=&e_{v\lambda(u)}+u-\lambda(u) &\qquad v\in\sigma,\ u\in \tau\ \textrm{and}\ \lambda(u)\notin\sigma\\
\end{aligned}
\end{equation}
where $e_{v\lambda(u)}$ are the standard basis vectors of $\ZZ^{\sigma\times(\widetilde\sigma\setminus\sigma)}$.
Then the closure of the image of $\Gr(k,\U)$ in the Pl\"ucker chart of $\bF_k(X_\A)$ corresponding to $\sigma$ is
the affine toric variety whose coordinate ring is
\[
\KK[S(\pi,\widetilde\sigma,\sigma)].
\]
In other words, on this Pl\"ucker chart, $Z_{\pi,k}$ is the above affine toric variety. Note that different choices of $\widetilde\sigma$ lead to isomorphic semigroups $S(\pi,\widetilde\sigma,\sigma)$.

If $\sigma$ is a face of $\A$, this chart of $Z_{\pi,k}$ has a (unique) torus fixed point, which is exactly the fixed point corresponding to $\sigma$ under Proposition \ref{prop:fixed}. We thus have
\begin{prop}\label{prop:fixed2}
The torus fixed points of $Z_{\pi,k}$ are in bijection with empty simplicial $k$-faces $\sigma$ of $\tau$ mapping to $k+1$ vertices of $\Delta_l$. 
\end{prop}

\begin{ex}\label{ex:birkhoff3}
Using the above local description of $Z_{\pi,k}$, we get explicit local descriptions for the irreducible componenets of $\bF_2(X_\A)$, where $\A$ is as in Examples \ref{ex:birkhoff1} and \ref{ex:birkhoff2}. Indeed, for Cayley structures of the form $\pi:\tau\to\Delta_3$, $\tau$ a facet of $\A$, $Z_\pi$ is just the single point $[L_\pi]$, and $Z_{\pi,2}$ is thus isomorphic to the Grassmannian $G(3,4)\cong\PP^3$.

 On the other hand, consider one of the maximal Cayley structures $\pi:\A\to\Delta_2$, for example, $\pi(u_i)=\pi(v_i)=e_i$. If we choose $\sigma=\{v_0,u_1,u_2\}$, the corresponding matrix $p_{vu}(t)$ has the form 
\[
\begin{blockarray}{ccccccc}
& v_0&u_1&u_2&u_0&v_1&v_2\\
\begin{block}{c(cccccc)}
v_0   &  1  &0  &0&t_1t_2 & 0 & 0\\
u_1   & 0&  1    &0& 0& t_1& 0 \\
u_2 &  0&0&  1 &0 & 0 & t_2\\
\end{block}
\end{blockarray}
\]
\end{ex}
\noindent for appropriate choice of parameters $t_1,t_2$. 
Hence, the component $Z_\pi$ is covered by copies of $\Aff^2$.

\begin{rem}\label{rem:chow}
A central computation of \cite{ilten:16a} involves understanding those $5$-planes of $\PP^{11}$ which are not contained in a coordinate hyperplane, yet are contained in the hypersurface $X$ on which
\[
f=x_0x_1x_2+x_3x_4x_5+x_6x_7x_8+x_9x_{10}x_{11}
\]
vanishes.
Setting $f_1=x_0x_1x_2+x_3x_4x_5$ and $f_2=x_6x_7x_8+x_9x_{10}x_{11}$, clearly if a $5$-plane $L$ is contained in the variety $V(f_1,f_2)=V(f_1)\cap V(f_2)$ on which $f_1$ and $f_2$ vanish, it is contained in $X$. But since $V(f_i)$ is just the cone over the variety $X_\A$ from Example \ref{ex:birkhoff3}, it is straightforward to see that any $5$-plane $L$ contained in $V(f_1,f_2)$ but not in a coordinate hyperplane comes in a $4$-dimensional family (up to certain permutations of the variables).

Somewhat surprisingly, after allowing for different choices of $f_1$ and $f_2$ by permuting the terms of $f$, all $5$-planes of $X$ not contained in a coordinate hyperplane arise in this fashion \cite[Proposition 3.2]{ilten:16a}. It would be interesting to provide a non-computational proof of this fact, and understand more generally in what situations linear spaces of a special (yet non-toric) hypersurface arise in this toric fashion.
\end{rem}

\begin{prop}\label{prop:contain}
Let $\pi:\tau \to \Delta_{l}$ and $\pi':\tau'\to \Delta_{l'}$ be Cayley structures for $l,l'\geq k\geq 1$. Then $Z_{\pi,k}$ contains $Z_{\pi',k}$ if and only if $\pi\mgeq \pi'$.
\end{prop}
\begin{proof}
It is straightforward to check from the construction of $Z_{\pi,l}$ that $\pi\mgeq \pi'$ implies $Z_{\pi',k}\subset Z_{\pi,k}$.
On the other hand, suppose that  $Z_{\pi',k} \subset Z_{\pi,k'}$. Since every $k$-plane of $L_{\pi'}$ must be contained in $X_{\tau}\subset X_\A$, we see that $\tau'\subset \tau$, so $\tau'$ is in fact a face of $\tau$.

To construct a map $\rho:\Delta_{l}\to\Delta_{l'}$ such that $(\rho\circ\pi)_{|\tau'}=\pi'$, define
$\rho(e_i)$ to be $e_0$ if $e_i\notin \pi(\tau')$. Otherwise, choose some $u\in \tau'$ such that $\pi(u)=e_i$ and set $\rho(e_i)=\pi'(u)$. We must show that for $u,u'\in \tau'$, if $\pi(u)=\pi(u')$, then $\pi'(u)=\pi'(u')$.

If $\pi'(u)\neq \pi'(u')$, note that for any matrix parametrizing a sufficiently general $k$-plane of $L_{\pi'}$, the two columns corresponding to $u$ and $u'$ will be linearly independent. On the other hand, consider any empty simplicial $k$-face $\sigma$ of $\tau'$ such that $\pi(\sigma)$ contains $\pi(u)$. By inspecting the  matrix $(q_{iu}(t))$ parametrizing $Z_{\pi,k}$  on the local chart corresponding to $\sigma$, we note that the two columns corresponding to $u$ and $u'$ will be linearly dependent. Hence, we must have $\pi'(u)=\pi'(u')$.
\end{proof}

\section{Proof of Main Result}\label{sec:proof}
Let $L$ be a $k$-dimensional linear subspace of $X_\A$. We will show that $[L]$ is a point of $Z_{\pi,k}$ for some Cayley structure $\pi$. Hence, every point of $\bF_k(X_\A)$ is contained in some $Z_{\pi,k}$, so Theorem \ref{thm:main} will then follow from Proposition \ref{prop:contain}.

We can represent $L$ as the rowspan of a full-rank $((k+1)\times \#\A)$-matrix $(\alpha_{iu})$ with $i=0,\ldots,k$ and columns indexed by $u\in \A$. We let 
\[
y_u=\sum_{i=0}^k\alpha_{iu}y_i
\]
be a linear form in indeterminates $y_0,\ldots,y_k$. The criterion that $L$ be contained in $X_\A$ is exactly the condition that for  any any relation
\[
\sum_{u\in\A} a_u u=\sum_{u\in\A} b_u u
\]
with $a_u,b_u\in \ZZ_{\geq 0}$ and $\sum a_u=\sum b_u$, we have
\begin{equation}\label{eqn:key}
\prod_{u\in\A} y_u^{a_u} =\prod_{u\in\A} y_u^{b_u}.
\end{equation}
Indeed, the defining equations of $X_\A$ are binomials corresponding to such affine relations above.
  
Let $\tau$ consist of those $u\in\A$ such that $y_u$ is non-zero.
\begin{lemma}
The set $\tau$ is a face of $\A$.
\end{lemma}
\begin{proof}
Considering a generic point $\eta$ of $L$ with coordinates $\eta_u$ for $u\in\A$, the set $\tau$ consists of exactly those $u$ with $\eta_u\neq 0$. On the other hand, $\eta$ is contained in a torus orbit $\co$ of $X_\A$. Such an orbit $\co$ correspond to a face $\widetilde\sigma$ of $\A$, see \S\ref{sec:toric}. The points $\zeta=(\zeta_u)$  contained in the orbit $\co$ are precisely those points for which  $\zeta_u\neq 0$ if and only if $u\in \widetilde\sigma$. The claim now follows.   
\end{proof}

Let $V$ be the vector space of linear forms in $\KK[y_0,\ldots,y_k]$. Then we have a map 
\begin{align*}
\rho:\tau&\to \PP(V)\\
u&\mapsto \langle{y_u}\rangle
\end{align*}
where $\langle{y_u}\rangle$ is the line in $\PP(V)$ spanned by $y_u$. Composing $\rho$ with a bijection between $\rho(\tau)$ and the elements of $\Delta_l$ for $l=\#\rho(\tau)-1$, we obtain a map
$\pi:\tau\to\Delta_l$. Note that since the rank of $(\alpha_{iu})$ is $k+1$, we must have $l\geq k$. 

\begin{lemma}\label{lemma:Cayley}
The map $\pi$ is a Cayley structure.
\end{lemma}
\begin{proof}
Consider any affine relation 
\[
\sum_{u\in\tau} a_u u=\sum_{u\in\tau} b_u u
\]
with $a_u,b_u\in \ZZ_{\geq 0}$ and $\sum a_u=\sum b_u$. Then by equation \eqref{eqn:key} and the fact that $\KK[y_0,\ldots,y_k]$ is a unique factorization domain, we must have 
\[
\sum_{\rho(u)=v} a_u =\sum_{\rho(u)=v} b_u 
\]
for each $v\in\rho(\tau)$. But this implies that
\[
\sum_{u\in\tau} a_u \pi(u)=\sum_{u\in\tau} b_u \pi(u).
\]
Hence, $\pi$ is a Cayley structure.
\end{proof}

\begin{rem}\label{rem:cover}
Lemma \ref{lemma:Cayley} provides a straightforward proof that $X_\A$ is covered by $k$-planes if and only if there is a Cayley structure $\pi:\A\to \Delta_k$. Indeed, if $X_\A$ is covered by $k$-planes, there is a $k$-plane $L$ containing a point in the dense torus orbit of $X_\A$. For this plane $L$, the set $\tau$ as above must be all of $\A$, and by the lemma we obtain a Cayley structure $\pi:\A\to\Delta_l$ for some $l\geq k$. We can compose this with any affine surjection $\Delta_l\to\Delta_k$ to get the desired Cayley structure. On the other hand, if $\pi:\A\to\Delta_k$ is a Cayley structure, the set of $k$-planes parametrized by $Z_{\pi}$ clearly covers $X_\A$. We thus recover \cite[Theorem 1.2]{ito:15a} for polarized normal toric varieties $(X,\mcL)$ in the special case that $\mcL$ is ample. 
\end{rem}

Returning to our proof of Theorem \ref{thm:main}, we claim now that the point $[L]\in\bF_k(X_\A)$ lies in $Z_{\pi,k}\subset\bF_k(X_\A)$. Choose $\sigma\subset \tau$ such that the $(k+1)\times (k+1)$ matrix $(\alpha_{iv})_{v\in \sigma}$ is invertible.
We can thus represent $L$ by a matrix $(\alpha_{vu})$ with rows and columns indexed by $v\in \sigma,u\in \A$ such that the square submatrix $(\alpha_{vu})_{v,u\in\sigma}$ is the identity matrix. Similarly,  the forms $y_v$, $v\in \sigma$ are a basis of the vector space $V$ and we have
$y_u=\sum_{v\in\sigma} \alpha_{vu}y_v$.

Consider any $\widetilde\sigma\subset \tau$ such that $\sigma\subset \widetilde\sigma$ and $\pi$ gives a bijection between $\widetilde\sigma$ and the vertices of $\Delta_l$.
As we saw in \S\ref{sec:comps}, the choice of $\sigma$  gives an affine chart $U$ of $Z_{\pi,k}$. 
This chart $U$ is the toric variety corresponding to the semigroup $S(\pi,\widetilde\sigma,\sigma)$ generated by $\gamma(u)$ and $\gamma(v,u)$ as in Equation \eqref{eqn:gens}.
Equivalently, this chart is the affine toric variety parametrized by the entries of the matrix from equation \eqref{eqn:matrix}.

We will show that $[L]$ is a point in this chart. 
 Now, all relations among these parametrizing functions are generated from those binomial equations $f=0$ corresponding to relations of the form
\begin{equation}\label{eqn:general}
\begin{aligned}
\sum_{\lambda(u)\in\sigma} a_u \gamma(u)
+\sum_{\substack{\lambda(u)\notin\sigma\\v\in\sigma}} a_u \gamma(v,u)
=\sum_{\lambda(u)\in\sigma} b_u \gamma(u)
+\sum_{\substack{\lambda(u)\notin\sigma\\v\in\sigma}} b_u \gamma(v,u)
\end{aligned}
\end{equation}
in $M\times\ZZ^{\sigma\times(\widetilde\sigma\setminus\sigma)}$.
Here, $a_u,b_u\in\ZZ_{\geq 0}$. We need to show that all such binomials $f$ vanish at the point $[L]=(\alpha_{vu})$, that is, that
\begin{equation}\label{eqn:firstsat}
\prod_{\lambda(u)\in\sigma} \alpha_{\lambda(u)u}^{a_u}\prod_{\substack{\lambda(u)\in\widetilde\sigma\setminus\sigma\\v\in\sigma}} \alpha_{vu}^{a_u}=
\prod_{\lambda(u)\in\sigma} \alpha_{\lambda(u)u}^{b_u}\prod_{\substack{\lambda(u)\in\widetilde\sigma\setminus\sigma\\v\in\sigma}} \alpha_{vu}^{b_u}.
\end{equation}
 
Fix an order $<$ on the elements of $\sigma$. We define a map $\kappa:\widetilde{\sigma}\setminus\sigma\to\sigma$ as follows: for any $v\in \widetilde\sigma\setminus\sigma$, let $\kappa(v)$ be the smallest element of $\sigma$ such that $\alpha_{\kappa(v)v}\neq 0$. 
For any $u$ with $\lambda(u)=v\in\widetilde\sigma\setminus\sigma$, we set $\kappa(u)=\kappa(\lambda(u))$.
Our conditions on $\kappa(u)$ imply that for any set of natural numbers $a_u$, the coefficient of 
\[
\prod_{\lambda(u)\in\widetilde\sigma\setminus\sigma}y_{\kappa(u)}^{a_u} \qquad
\mathrm{in}
\prod_{\lambda(u)\in\widetilde\sigma\setminus\sigma}y_{u}^{a_u}
\]
is
\begin{equation}\label{eqn:monomial}
\prod_{\lambda(u)\in\widetilde\sigma\setminus\sigma}\alpha_{\kappa(u)u}^{a_u}.
\end{equation}

Some special affine relations are those of the form
\begin{equation*}
\gamma(v,u)+\gamma(\kappa(u), \lambda(u))=
\gamma(\kappa(u),u)+\gamma(v,\lambda(u))
\end{equation*}
for $v\in\sigma$ and $u\in\tau$ with $\lambda(u)\notin\sigma$.
The corresponding binomial evaluated at $[L]$ is exactly
\[
\alpha_{vu}\alpha_{\kappa(u)\lambda(u)}-\alpha_{\kappa(u) u}\alpha_{v\lambda(u)}.
\]
This expression is indeed zero: $\pi(u)=\pi(\lambda(u))$ implies that $y_u$ and $y_{\lambda(u)}$ are linearly dependent.

Since $\alpha_{\kappa(v) \lambda(u)}\neq 0$, we can rewrite
\[
\alpha_{vu}=\frac{\alpha_{\kappa(u) u}\alpha_{v\lambda(u)}}{\alpha_{\kappa(u)\lambda(u)}}.
\]
Thus, it suffices to consider relations \eqref{eqn:general} where the only $\gamma(v,u)$ terms appearing are those with either $v=\kappa(u)$ or $u=\lambda(u)$.
Indeed, if \eqref{eqn:firstsat} is satisfied for such relations, we obtain equality for arbitrary relations by repeated substitutions using the above expressions for $\alpha_{vu}$.  
Furthermore, for such a relation, the only terms containing an $e_{v\lambda(u)}$ component for $v\neq \kappa(u)$ are $\gamma(v,\lambda(u))$. We can cancel these from both sides of the relation to assume that the only $\gamma(v,u)$ terms appearing are those with  $v=\kappa(u)$.

We can thus reduce our general relations \eqref{eqn:general} to 
\begin{equation}\label{eqn:special}
\begin{aligned}
\sum_{\lambda(u)\in\sigma} a_u \gamma(u)
+\sum_{\substack{\lambda(u)\notin\sigma}} a_u \gamma(\kappa(u),u)
=\sum_{\lambda(u)\in\sigma} b_u \gamma(u)
+\sum_{\substack{\lambda(u)\notin\sigma}} b_u \gamma(\kappa(u),u).
\end{aligned}
\end{equation}
We must then show that
\begin{equation}\label{eqn:want}
\prod_{\lambda(u)\in\sigma} \alpha_{\lambda(u)u}^{a_u}\prod_{\lambda(u)\in\widetilde\sigma\setminus\sigma} \alpha_{\kappa(u)u}^{a_u}=
\prod_{\lambda(u)\in\sigma} \alpha_{\lambda(u)u}^{b_u}\prod_{\lambda(u)\in\widetilde\sigma\setminus\sigma} \alpha_{\kappa(u)u}^{b_u}.
\end{equation}

Projecting \eqref{eqn:special} to $\ZZ^{\sigma\times(\widetilde\sigma\setminus\sigma)}$, we obtain that for any $v\in\widetilde\sigma\setminus\sigma$,
\[
\sum_{\lambda(u)=v} a_u=\sum_{\lambda(u)=v} b_u.
\]
Indeed, for such $v$, the coefficient of $e_{\kappa(v)v}$ in the left hand side of the projection of \eqref{eqn:special} is $\sum_{\lambda(u)=v} a_u$, and a similar claim holds for the right hand side. 
Thus,
\[
\prod_{\lambda(u)\notin\sigma} \alpha_{\kappa(u)\lambda(u)}^{a_u}=\prod_{\lambda(u)\notin\sigma} \alpha_{\kappa(u)\lambda(u)}^{b_u}.
\]
We denote this quantity by $c$; note that by our construction of $\kappa$,  $c\neq 0$.

On the other hand, rearranging and projecting \eqref{eqn:special} to $M$, we have the affine relation
\begin{align*}
\sum_{u\in\tau} a_u u+b_u \lambda(u)
=\sum_{u\in\tau} b_u u+a_u\lambda(u)
\end{align*}
among the elements of $\tau\subset \A$. Passing to the induced relation on the $y_u$ from equation \eqref{eqn:key} gives 
\begin{equation*}
\prod_{u\in\tau} y_u^{a_u}\prod_{u\in\tau} y_{\lambda(u)}^{b_u}=\prod_{u\in\tau} y_{\lambda(u)}^{a_u} \prod_{u\in\tau} y_u^{b_u}.
\end{equation*}
We view both sides of this equation as polynomials in the $y_v$, $v\in\sigma$. 
Consider the coefficients on both sides of this equation for the monomial
\[
\prod_{v\in\sigma}y_v^{d_v}\cdot\prod_{v\in\widetilde\sigma\setminus\sigma} y_{\kappa(v)}^{d_v}
\]
where 
\begin{align*}
d_v=\sum_{\pi(u)=\pi(v)} a_u+b_u.
\end{align*}
But by \eqref{eqn:monomial} these are just $c$ times the left and right hand sides of Equation \eqref{eqn:want}. We conclude that $[L]$ satisfies the necessary binomial equations and is thus contained in $Z_{\pi,k}$. 

We have now seen that $\bF_k(X_\A)=\bigcup Z_{\pi,k}$, where the union is taken over all maximal Cayley structures $\pi$. But then each $Z_{\pi,k}$ must be an irreducible component by Proposition \ref{prop:contain}.
This completes our proof of Theorem \ref{thm:main}.

\section{Dimension and Global Description}
Having now understood the component structure of $\bF_k(X_\A)$, we begin to study some properties of the components. We begin with a description of the dimension of the varieties $Z_{\pi,k}$.

\begin{prop}\label{prop:dim}
Let $\tau$ be an $m$-face of $\A$ and $\pi:\tau\to\Delta_l$ a Cayley structure for some $l\geq k$. Then the subvariety $Z_{\pi,k}$ of $\bF_k(X_\A)$ has dimension $m-l+(k+1)(l-k)$.
\end{prop}
\begin{proof}
Let $\widetilde \sigma$ be any $l$-face of $\tau$ surjecting onto $\Delta_l$, and $\sigma$ any $k$-face of $\widetilde \sigma$. By the local description of $Z_{\pi,k}$ from \ref{sec:comps}, its dimension is given by the dimension of the semigroup $S(\pi,\widetilde\sigma,\sigma)$ contained in $M\times \ZZ^{\sigma\times(\widetilde\sigma\setminus\sigma)}$. The projection $\overline S$ of this semigroup to the $M$ factor is generated by $u-\lambda(u)$ for $u\in \tau$.  The dimension of $\overline S$ is $m-l$, since $\pi$ is an affine linear map, which when extended to a map from the affine span of $\tau$ to $\RR^{l+1}$ has $(m-l)$-dimensional fibers.

We thus see that the dimension of $S(\pi,\widetilde\sigma,\sigma)$ is at most 
\[m-l+\#(\sigma\times(\widetilde\sigma\setminus\sigma)=
m-l+(k+1)(l-k).
\]
But this dimension is actually achieved, since in the projection $S(\pi,\widetilde\sigma,\sigma)\to\overline S$ the fiber over $0$ consists of the $(k+1)(l-k)$ linearly independent $e_{vu}$ for $v\in\sigma$, $u\in\widetilde\sigma\setminus\sigma$.
\end{proof}

For a Cayley structure $\pi:\tau\to\Delta_l$, $Z_\pi$ is the $T$-orbit closure of the point $[L_\pi]$ in the Fano scheme $\bF_l(X_\A)$. As such, it is a (potentially non-normal) toric variety. Embedding in $\bF_l(X_\A)$ in projective space by the Pl\"ucker embedding we arrive at an equivariant embedding of $Z_\pi$. This can be described globally as follows. Let $\A_\pi\subset M$ be
\[
\A_\pi=\{u_0+\ldots+u_l\in M\ |\ u_i\in\pi^{-1}(e_i)\}
\]
\begin{thm}\label{thm:polytope}
In its Pl\"ucker embedding, $Z_\pi$ is the toric variety $X_{A_\pi}$. 
\end{thm}
\begin{proof}
We simply consider the action of $T$ on each Pl\"ucker coordinate. As noted in \S\ref{sec:comps}, the non-zero Pl\"ucker coordinates for $Z_\pi$ correspond to choosing one element from each fiber $\pi^{-1}(e_i)$, that is, exactly the tuples $(u_0,\ldots,u_l)$ appearing in the construction of $\A_\pi$. Taking into account the action of $T$ on $L_\pi$ as described in \S\ref{sec:comps}, we see that $T$ acts on the coordinate corresponding to the tuple $(u_0,\ldots,u_l)$ with weight $u_0+\ldots+u_l$.
\end{proof}

\begin{rem}
Note that the action of $T$ on $Z_\pi$ is not faithful: a faithful action is given by the torus $\overline T$ whose characters are generated by the kernel of $\pi$. However, the canonical embedding of $Z_\pi$ in projective space does not possess a canonical $\overline T$-equivariant structure. 
\end{rem}

\begin{rem}
Although the components $Z_{\pi,k}$ are locally toric, they are in general not globally toric. 
\end{rem}

\begin{ex}
Continuing Example \ref{ex:birkhoff3}, we see had already seen that $\bF_3(X_\A)$ has $6$ $2$-dimensional components, and $9$ $3$-dimensional components, which agrees with the contents of Proposition \ref{prop:dim}. As noted before, the $3$-dimensional components are all isomorphic to $\PP^3$. We may apply Theorem \ref{thm:polytope} to determine that for one of the $2$-dimensional components, $\A_\pi$ is given by
\begin{align*}
&u_0+u_1+u_2\qquad &u_0+u_1+v_2 \qquad &u_0+v_1+u_2\qquad &u_0+v_1+v_2&\\
&v_0+u_1+u_2\qquad &v_0+u_1+v_2 \qquad &v_0+v_1+u_2\qquad &v_0+v_1+v_2&
\end{align*}
up to some permutation of the $v_j$. This configuration of lattice points is isomorphic to the set of lattice points in $\ZZ^2$ given by the columns of 
\[
\left(\begin{array}{c c c c c c c c}
0 & 0 & 0 & 1 & 1 & -1 & -1\\
0 & 1 & -1& 0 & 1 & 0 & -1
\end{array}\right).
\]
The corresponding toric variety is just $\Bl_3\PP^2$, the blowup of $\PP^2$ in $3$ points, embedded by its anticanonical divisor.
\end{ex}

\section{Smoothness of $\bF_k(X_\A)$ }
\begin{defn}\label{defn:smooth}
Let $\sigma$ be an empty simplicial $k$-face of $\A$. We say that $\A$ is \emph{smooth at $\sigma$} if the semigroup generated by $\A-\sigma$ is isomorphic to $\ZZ^{k}\times\ZZ_{\geq0}^{n-k}$, where $n$ is the dimension of $\A$.
\end{defn}

\begin{rem}\label{rem:regular}
The set $\A$ is smooth at $\sigma$ if and only if $X_\sigma\subset X_\A$ is not contained in the singular locus of $X$. Indeed, $X_\sigma$ is the orbit closure of the distinguished point $\eta$ (cf. \cite[\S3.2]{CLS}) of the affine patch $U=\spec \KK [\NN\cdot (\A-\sigma)]$, so it suffices to check if the point $\eta$ is in the singular locus of $X_\A$. Furthermore, $U$ is the smallest torus invariant open set containing $\eta$. Since the singular locus of $X_\A$ is torus invariant, $X_\A$ is smooth at $\eta$ if and only if $U$ is smooth. This is equivalent to $\NN\cdot(\A-\sigma)$ being isomorphic to  $\ZZ^{k}\times\ZZ_{\geq0}^{n-k}$, see \cite[\S1.3]{CLS}.
\end{rem}

\begin{thm}\label{thm:regular}
Let $\sigma$ be an empty simplicial $k$-face of $\A$. Suppose that $\A$ is smooth at $\sigma$. Then any component $Z$ of $\bF_k(X_\A)$ containing the point $[L_\sigma]$ corresponding to $\sigma$ is smooth at $[L_\sigma]$.
\end{thm}
\begin{proof}
We know by Theorem \ref{thm:main} and Proposition \ref{prop:fixed2} that the component $Z$ is of the form $Z_{\pi,k}$ for some Cayley structure $\pi:\tau\to \Delta_{l}$, $l\geq k$ with $\sigma$ a face of $\tau$ mapping to $k$ vertices of $\Delta_l$.
Let $m$ denote the dimension of $\tau$. Note that the set of all $u\in\tau$ with $\lambda(u)\in\sigma$ is a $m+k-l$-dimensional face $\omega$ of $\tau$. 
Since $\A$ is smooth at $\sigma$, we find $u_1,\ldots,u_{n-k}\in\A$ such that every $u\in\A$ can be written uniquely as 
\begin{equation}\label{eqn:u}
u=\sum_i a_i u_i+\sum_{v\in\sigma} b_vv
\end{equation}
with $a_i\in\ZZ_{\geq0}$, $b_v\in \ZZ$, and $\sum_i a_i+\sum_v b_v=1$. 
After reordering the $u_i$, we can assume that exactly $u_1,\ldots,u_{m-k}$ are in $\tau$, and $u_1,\ldots,u_{m-l}$ are in $\omega$. This follows from the fact that both sets are faces of $\A$.
Furthermore, $u$ is in $\tau$ (respectively $\omega$) if and only if the representation \eqref{eqn:u}, has $a_i=0$ for $i>m-k$ (respectively $i>m-l$).

To study the Pl\"ucker chart of $Z_{\pi,k}$ containing $[L_\sigma]$, we consider the set $\widetilde\sigma=\sigma\cup\{u_{m-l+1},\ldots,u_{m-k}\}$. This has exactly $l+1$ elements and surjects (via $\pi$) onto the vertices of $\Delta_l$. Indeed, consider any $u\in\tau$ with $\pi(u)=e_j\notin\pi(\sigma)$. Then since $\pi$ is affine, using \eqref{eqn:u} we have
\[
e_j=\sum_{i\leq m-k} a_i\pi(u_i)+\sum_{v\in\sigma}b_v\pi(v).
\]
The only way this can happen is if $\pi(u_i)=e_j$ for some $i$, since $a_i\geq 0$.
Thus, as in \S\ref{sec:comps}, the Pl\"ucker chart of $Z_{\pi,k}$ containing $x_\sigma$ is the toric variety associated to the semigroup $S(\pi,\widetilde\sigma,\sigma)$.

We will show that $S(\pi,\widetilde\sigma,\sigma)$ is generated by $\gamma(u_i)=u_i-\lambda(u_i)$ for $i<m-l$, and $e_{vw}$ for $w\in \widetilde\sigma\setminus\sigma$, $v\in\sigma$. If this claim is true, we are done, since these generators are linearly independent; this implies the smoothness of $U$.

Since \eqref{eqn:u} is an affine relation, we have that
\begin{equation*}
\pi(u)=\sum_i a_i \pi(u_i)+\sum_{v\in\sigma} b_v\pi(v)
\end{equation*}
for any $u\in\tau$.
But the relations among the points of $\sigma$ are exactly the same as those among the points of $\Delta_l$, hence
\begin{equation*}
\lambda(u)=\sum_i a_i \lambda(u_i)+\sum_{v\in\sigma} b_vv
\end{equation*}
and we obtain
\[
u-\lambda(u)=\sum_{i\leq l-k} a_i(u_i-\lambda(u_i))+\sum_{i > l-k} a_i(u_i-u_i)=\sum_{i\leq l-k} a_i(u_i-\lambda(u_i)).
\]
In particular, if $u\in \omega$, then $\gamma(u)$ is generated by the $\gamma(u_i)$, $i=1,\ldots,m-l$.
Furthermore, if $u\in\tau\setminus\omega$, then $\gamma(v,u)$ is generated by $e_{v\lambda(u)}$ along with the $\gamma(u_i)$, $i=1,\ldots,m-l$. 

To complete the proof, we note that our desired generators are indeed in $S(\pi,\widetilde\sigma,\sigma)$: $\gamma(u_i)$ obviously so, and $e_{v\lambda(u)}$ since $\gamma(v,\lambda(u))=e_{v\lambda(u)}$.
\end{proof}

\begin{cor}\label{cor:regular}
Suppose that the dimension of the singular locus of $X_\A$ is at most $k-1$. Then in its reduced structure, each component of $\bF_k(X_\A)$ is smooth.
\end{cor}
\begin{proof}
A component $Z_{\pi,k}$ is smooth if and only if it is smooth at its toric fixed points. The claim now follows from Theorem \ref{thm:regular} and Remark \ref{rem:regular}.
\end{proof}

\section{Intersections and Connectedness}
We can also completely describe the intersection behaviour of the subvarieties $Z_{\pi,k}$ of $\bF_k(X_\A)$. 
\begin{thm}\label{thm:intersect}
Let $\pi_1:\tau_1\to\Delta_{l_1}$, $\pi_2:\tau_\to\Delta_{l_2}$ be Cayley structures with $l,l'\geq k$. Then 
\[
Z_{\pi_1,k}\cap Z_{\pi_2,k}=\bigcup_{\substack{\pi\preceq \pi_i\\i=1,2}} Z_{\pi,k}.
\]
\end{thm}
\begin{proof}
By Proposition \ref{prop:contain}, the right hand side is clearly contained in the left.
On the other hand, consider some linear space $L$ such that $[L]\in Z_{\pi_1,k}\cap Z_{\pi_2,k}$.
We may construct a Cayley structure $\pi:\tau\to\Delta_l$ as in \S\ref{sec:proof} such that $[L]\in Z_{\pi,k}$. Using notation as in \S\ref{sec:proof}, $u\in\tau$ only if $y_u\neq 0$, but since $[L]\in Z_{\pi_1,k}\cap Z_{\pi_2,k}$ this can only happen if $\tau\subset \tau_1\cap\tau_2$. Hence, $\tau$ is a face of $\tau_1$ and $\tau_2$. Furthermore, $\pi(u)\neq \pi(v)$ if and only if $y_u$ and $y_v$ are not linearly dependent. But again since $[L]$ is contained in the left hand side above, we must also have $\pi_i(u)\neq \pi_i(v)$. Hence, we can find maps $\rho_i:\Delta_{l_i}\to\Delta_l$ such that $(\rho_i\circ\pi_i)_{\tau}=\pi$, so $Z_{\pi,k}$ contains $[L]$ and appears on the right hand side above.
\end{proof}

To describe the connected components of $\bF_k(\A)$, we construct a graph $\Gamma_\A$ based on the combinatorics of $\A$, with the property that its connected components are in bijection with those of $\bF_k(\A)$.
The vertex set of $\Gamma_\A$ consists of the set of maximal Cayley structures $\pi:\tau\to\Delta_{l}$ for $\tau\prec \A$ and $l\geq k$. We connect vertices $\pi:\tau\to \Delta_l$ and $\pi':\tau'\to\Delta_{l'}$ by an edge in $\Gamma_\A$ if there exists an empty simplicial $k$-face $\sigma\prec\tau\cap\tau'$ such that $\pi$ and $\pi'$ are injective on $\sigma$.

\begin{thm}\label{thm:connected}
Irreducible components of $\bF_k(X_\A)$ are in bijection with vertices of $\Gamma_\A$. Two components intersect if and only if the corresponding vertices are connected by an edge. Two components are in the same connected component if and only if the corresponding vertices are in the same connected component of $\Gamma_\A$. In particular, the number of connected components of $\bF_k(X_\A)$ equals the number of connected components of $\Gamma_\A$.
\end{thm}

\begin{proof}
The first claim is just a restatement of Theorem \ref{thm:main}. For the second claim, note that if two components intersect, the intersection is $T$-invariant and projective, and thus must contain a torus fixed point. The claim now follows from Proposition \ref{prop:fixed2}. Two components $Z,Z'$ are in the same connected component if and only if there is a sequence of components $Z=Z_0,Z_1,\ldots,Z_j=Z'$ with $Z_i$ and $Z_{i+1}$ intersecting; the remaining claims follow.
\end{proof}

\begin{ex}
Continuing Example \ref{ex:birkhoff3}, it is easy to check via Theorem \ref{thm:connected} that $\bF_2(X_\A)$ is connected. Perhaps more interestingly, we may observe how the irreducible components intersect using Theorem \ref{thm:intersect}. Each $\PP^3$ component intersects $4$ other $\PP^3$ components in a point, missing the other $4$. Each $\PP^3$ component intersects $4$ $\Bl_3 \PP^2$ components in a $\PP^1$, and misses the other $2$.
The $\Bl_3 \PP^2$ components do not pairwise intersect.
\end{ex}

\section{Scheme Structure and Multiplicities}\label{sec:mult}
In this section, we will completely describe the scheme structure of $\bF_k(X_\A)$ locally around toric fixed points in the special case that $k=\dim X_\A-1$ and $X_\A$ is smooth in codimension one.

We now fix notation for the rest of this section.
Assume that $k=\dim X_\A-1$. Recall from Proposition \ref{prop:fixed} that the fixed points of $\bF_k(X_\A)$ correspond to faces of $\A$ which are empty $k$-simplices. Let us fix such a facet $\sigma=\{v_0,\ldots,v_k\}$ of $\A$.
Assume that $\A$ is smooth at $\sigma$ (cf. Definition \ref{defn:smooth}). This is certainly satisfied if $X_\A$ is smooth in codimension one. 
We will proceed to describe the coordinate ring $R_\sigma$ of the Pl\"ucker chart containing the fixed point corresponding to $\sigma$.

By the smoothness assumption, there exists some $w\in\A$ such that every $u\in \A$ can be written uniquely as 
\begin{equation}\label{eqn:c}
u=v_0+h\cdot (w-v_0)-\sum_{i=1}^k c_i (v_i-v_0)
\end{equation}
for $h\in\ZZ_{\geq 0}$, $c_i\in \ZZ$. We fix such an element $w$, and call $h$ the \emph{height} of $u$.
We then set
\[c_0=h-1-\sum_{i>0} c_i.\]
To stress the dependence of $c_i$ and $h$ on $u$, we will sometimes write $c_i(u)$ and $h(u)$. Note that $c_i(\cdot)$ and $h(\cdot)$ preserve affine relations among elements of $\A$.

Let $R=\KK[z_0,\ldots,z_k]$. For each element $u\in \A\setminus (\sigma\cup\{w\})$, we will describe a set of elements $S_u$ in $\ZZ_{\geq 0}^{k+1}$ such that a basis for $R_\sigma$ is given by $z^\alpha$ for  $\alpha\in \bigcap S_u$. Let $e_0\ldots,e_k$ be the standard basis of $\ZZ_{\geq 0}^{k+1}$. 
 For $\alpha=(\alpha_0,\ldots,\alpha_k)\in \ZZ_{\geq 0}^{k+1}$, set $|\alpha|=\sum_i \alpha_i$. 
Define $S_{<h}$ by $S_{<h}=\{\alpha\ ;\ |\alpha|<h\}$.

\begin{figure}
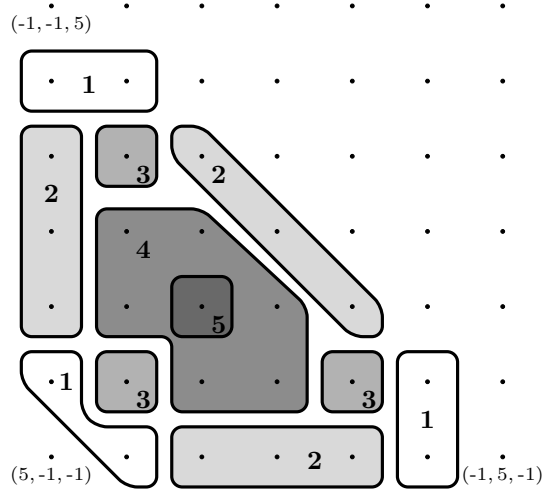

\casepic
\caption{Cases of Definition \ref{defn:S}.}
\label{fig:cases}
\end{figure}

\begin{defn}\label{defn:S}
Fix $u\in \A\setminus (\sigma\cup\{w\})$ and let $h,c_i$ be as in \eqref{eqn:c}.  We define $S_u$ as follows:
\begin{enumerate}
\item If for some $j\neq l$, $c_i=0$ for all $i\neq j,l$ and  $c_j=-1$ then \[S_u=\{\alpha\ |\ \alpha_i=0\ (i\neq l)\}\cup S_{<h};\] \label{case:simplex}
\item If for some $j$, $0\leq c_i<h$ for all $i\neq j$ and $c_j=-1$, then \[S_u=\{c+e_j \}\cup S_{<h};\]
\item If for some $j$, $c_i=0$ for all $i\neq j$, then \[S_u=\{\alpha\ |\ \sum_{i\neq j} \alpha_i \leq 1\}\cup S_{<h};\] 
\item If for some $j\neq l$, $c_i=0$ for all $i\neq j,l$ and $c_j,c_l>0$, then
\[S_u= \{c+e_0,\ldots,c+e_k,c+e_j+e_l\}\cup S_{<h};\]\label{case:middle} 
\item If $c_i\geq 0$ for all $i\neq j,l$ with at least three non-zero, then
\[S_u= \{c+e_0,\ldots,c+e_k\}\cup S_{<h};\] 
\item Otherwise, $S_u=S_{<h}$.\label{case:other}
\end{enumerate}
See Figure \ref{fig:cases} for the hyperplane slice $h(u)=4$ illustrating in coordinates $c_0,c_1,c_2$ the possible cases when $k=2$.
\end{defn}

We can now completely describe the local structure of the scheme $\bF_k(X_\A)$ in combinatorial terms.
\begin{thm}\label{thm:scheme}
With notation as above, a basis for the coordinate ring $R_\sigma$ is given by $z^\alpha$ for  $\alpha\in \bigcap_{u\in\A\setminus\sigma} S_u$. Multiplication is defined by $z^\alpha\cdot z^\beta=z^{\alpha+\beta}$, where we set $z^{\alpha+\beta}=0$ if $\alpha+\beta\notin \bigcap S_u$.
\end{thm}

In order to prove the theorem, we first set up some further notation. Let $y_0,\ldots,y_k$ be indeterminates, and for any $u\in\A\setminus\sigma$ and any $i=0,\ldots k$, let $z_{iu}$ be an indeterminate. We set $y_u=\sum z_{iu}y_i$, where if $u=v_j\in \sigma$, $z_{iu}$ is $0$ for $i\neq j$ and $1$ for $i=j$. Locally, the coordinate ring of $\bF_k(X_\A)$ may be expressed as a quotient of $\KK[z_{iu}]$ by the ideal obtained from the conditions on the $z_{iu}$ necessary to satisfy
\begin{equation}\label{eqn:rel}
\prod y_u^{a_u}=\prod y_u^{b_u}.
\end{equation}
for any affine relation $\sum a_uu=\sum b_uu$ on the elements of $\A$.
For any $i=0,\ldots,k$, set $z_i=z_{iw}$.

\begin{lemma}\label{lemma:mon}
The monomials $y^{\alpha}$ appearing in  
\[y_u \prod_{i\geq 0}y_i^{c_i}\]
(viewed as a polynomial in the $y_i$ with coefficients $z_{ju}$)
are exactly those for which $\alpha\in S_u$ and $|\alpha|=h$.
\end{lemma}
\begin{proof}
The above expression is homogeneous of degree $h$, so clearly we must have $|\alpha|=h$. The lemma may be easily verified from Definition \ref{defn:S} case by case.
\end{proof}

\begin{lemma}\label{lemma:span}
The monomials $z^\alpha$ for  $\alpha\in \bigcap_{u\in\A\setminus\sigma} S_u$ span the ring $R_\sigma$ as a $\KK$-vector space.
\end{lemma}
\begin{proof}
Considering any $u\in\A$, the affine relation derived from \eqref{eqn:c} imposes the condition
\[
y_u=y_w^h \prod_{i\geq 0}y_i^{-c_i}.
\]
Hence, each $z_{iu}$ can be expressed linearly in terms of the $z_j$. Furthermore, by Lemma \ref{lemma:mon}, this same condition imposes that $z^\alpha=0$ if $\alpha\notin S_u$. Indeed, for $|\alpha|\leq h$ this follows directly from the lemma. For $|\alpha|>h$, one may verify case by case that each $S_u$ is the set of exponent vectors for the standard monomials of a monomial ideal generated in degree $h$, and the claim follows.
\end{proof}

\begin{proof}[Proof of Theorem \ref{thm:scheme}]
From Lemma \ref{lemma:span}, we see that $R_\sigma$ is a quotient of the $\KK$-algebra $R_\sigma'$ with basis $z^\alpha$ for $\alpha\in \bigcap_{u\in\A\setminus\sigma} S_u$, with multiplication as defined in the statement of the theorem. To show the equality of $R_\sigma$ and $R_\sigma'$, it suffices to verify the following:
For any affine relation $\sum a_uu=\sum b_uu$, the equality \eqref{eqn:rel} is satisfied in $R_\sigma'$, where we have written each $z_{iu}$ as a $\KK$-linear combination of the $z_i$ as in Lemma \ref{lemma:span}.

Now, our construction of $S_u$ ensures that 
\[
y_u=y_w^{h(u)} \prod_{i\geq 0}y_i^{-c_i(u)}
\]
holds in $R_\sigma'$, see the proof of Lemma \ref{lemma:span}.
But then for the above affine relation, we have
\begin{align*}
\prod y_u^{a_u}&=\prod_u\left( y_w^{a_uh(u)} \prod_{i\geq 0}y_i^{-a_uc_i(u)}\right)=
y_w^{\sum_u a_uh(u)}\prod_{i\geq 0} y_i^{-\sum_u a_uc_i(u)}\\
\\
&=y_w^{\sum_u b_uh(u)}\prod_{i\geq 0} y_i^{-\sum_u b_uc_i(u)}
=\prod y_u^{b_u}
\end{align*}
in $R_\sigma'$, exactly as desired.
\end{proof}

\begin{cor}
Assume that the fixed point corresponding to $\sigma$ is an isolated point of $\bF_k(X_\A)$. Then the corresponding $k$-plane $L_\sigma$ of $X_\A$ has multiplicity equal to 
\[
\#\bigcap_u S_u.
\]
\end{cor}
\begin{proof}
The multiplicity of $L_\sigma$ is just the degree of the isolated point $[L_\sigma]$ of $\bF_k(X_\A)$. But this equals the dimension of $R_\sigma$ as a $\KK$-vector space, and the result follows from Theorem \ref{thm:scheme}.
\end{proof}

\begin{cor}\label{cor:two}
Assume that the fixed point of $\bF_k(X_\A)$ corresponding to $\sigma$ is an isolated point, and assume that $\A$ contains some $w'\neq w$ in height one. Then the multiplicity of $L_\sigma$ in $X_\A$ is equal to the smallest natural number $m$ for which 
the set 
\[
\{ u\in\A\ |\ \textrm{height of}\ u\leq m\}
\]
is contained in $\sigma+\NN\cdot (w-v_i)$ for some $i$.
\end{cor}
\begin{proof}
Let $w'\neq w$ be another element of $\A$ in height one. If $w'$ is not contained in $\sigma+(w-v_i)$ for some $i$, then for $u=w'$ we are in case \ref{case:other} of Definition \ref{defn:S}, and it follows that $\#\bigcap S_u=1$, as desired.

If instead $w'$ is contained in $\sigma+(w-v_i)$ for some $i$, then for $u=w'$ we are in case \ref{case:simplex}, and so $S_{w'}$ consists of those $\alpha$ of the form $\alpha=\lambda\cdot e_i$. Hence, the multiplicity of $L_\sigma$ is the smallest number $m$ for which $m e_i\notin \bigcap S_u$. But this is easily seen to be the $m$ from above.
\end{proof}

\begin{rem}
Our Theorem \ref{thm:scheme} implies the results of \cite[Theorem 2.2]{ilten:14a}, in which the scheme structure of $\bF_1(X_\A)$ for projectively normal toric surfaces $X_\A$ is completely determined. Indeed, if there are at least two elements in height one over $\sigma$, then the argument from the proof of Corollary \ref{cor:two} applies and the result is straightforward.
If there is only a single lattice point in height one and $\A$ is the set of lattice points of a polytope $P$, a straightforward argument from \cite{ilten:14a} shows that either $\A$ is contained in a prism over $\sigma$, $\A$ contains a point in height two, or $\A$ contains a point in height three satisfying case \ref{case:middle} of Definition \ref{defn:S}. An analysis in each case using Theorem \ref{thm:scheme} yields the result.
\end{rem}

As we will see in the following example, the Fano scheme $\bF_k(X_\A)$ may have embedded components.
\begin{ex}
Consider the set $\A$ from Example \ref{ex:surface1} and Figure \ref{fig:ex}.
Let $\sigma=\{(0,0),(1,0)\}$. Then $\A$ is smooth at $\sigma$, so we may apply Theorem \ref{thm:scheme} to obtain the scheme structure of $\bF_k(X_\A)$ around $[L_\sigma]$.

Here, $w=(0,1)$, so the only other point we have is $u=(1,2)$. This $u$ is in case \ref{case:simplex} of Definition \ref{defn:S}, and $S_u=\{(1,0)\}\cup \{\lambda\cdot (0,1)\}_{\lambda\in\ZZ_{\geq 0}}$.
Hence, $\bF_k(X_\A)$ is locally an affine line with a fat point at the origin.
\end{ex}

\begin{rem}
While it seems rather challenging to give a palatable combinatorial description of the scheme structure of $\bF_k(X_\A)$ for $k<\dim X_\A-1$, a good first step would be to describe exactly what the embedded components of $\bF_k(X_\A)$ are. Conjecturally, these are all of the form $Z_{\pi,k}$ where $\pi$ is some non-maximal Cayley structure. 
\end{rem}

\bibliographystyle{amsalpha}
\renewcommand{\MR}[1]{\relax}
\bibliography{toric-fano}

\end{document}

%% file: figures.tex
\newcommand{\cspic}{\begin{tikzpicture}[scale=1]
\foreach \x in {-1,...,2}{
  \foreach \y in {0,...,2}{
    \node[draw,circle,inner sep=0.5pt,fill] at (\x,\y) {};
  }
}

\node[draw,circle,inner sep=2pt] at (0,0) {};
\node[draw,circle,inner sep=2pt] at (1,0) {};
\node[draw,circle,inner sep=2pt] at (0,1) {};
\node[draw,circle,inner sep=2pt] at (1,2) {};

\node at (0,-1) {$e_0$};
\node at (1,-1) {$e_1$};
\draw [->] (0,-.4) -- (0,-.7);
\draw [->] (1,-.4) -- (1,-.7);
\draw [rounded corners,line width=1] (-0.3,-0.3) -- (-0.3,1.3) -- (0.3,1.3) -- (0.3,-0.3) -- cycle;
\draw [rounded corners,line width=1] (0.7,-0.3) -- (0.7,.5) -- (1.15,.5) -- (1.15,1.5) -- (0.7,1.5) -- (0.7,2.3) -- (1.3,2.3) -- (1.3,-0.3) -- cycle;

\end{tikzpicture}}

\newcommand{\lppic}{\begin{tikzpicture}[scale=1]
\foreach \x in {-1,...,2}{
  \foreach \y in {-1,...,3}{
    \node[draw,circle,inner sep=0.5pt,fill] at (\x,\y) {};
  }
}

\node[draw,circle,inner sep=2pt] at (0,0) {};
\node[draw,circle,inner sep=2pt] at (1,0) {};
\node[draw,circle,inner sep=2pt] at (0,1) {};
\node[draw,circle,inner sep=2pt] at (1,2) {};

\end{tikzpicture}}

\newcommand{\casepic}{\begin{tikzpicture}[scale=1]

\draw [rounded corners,line width=1.2] (-2.4,-1) -- (-2.4,-0.6) -- (-1.6,-0.6) -- (-1.6,-1.4) -- (-1.4,-1.6) -- (-0.6,-1.6) -- (-0.6,-2.4) -- (-1,-2.4) -- cycle;
\draw [rounded corners,line width=1.2] (-2.4,2.6) -- (-2.4,3.4) -- (-0.6,3.4) -- (-0.6,2.6) -- cycle;
\draw [rounded corners,line width=1.2] (2.6,-2.4) -- (2.6,-0.6) -- (3.4,-0.6) -- (3.4,-2.4) -- cycle;
\node at (-1.8,-1) {\textbf{1}};
\node at (-1.5,2.95) {\textbf{1}};
\node at (3,-1.5) {\textbf{1}};

\draw [rounded corners,line width=1.2,fill opacity=0.25,fill=white!40!black] (-0.4,-2.4) -- (-0.4,-1.6) -- (2.4,-1.6) -- (2.4,-2.4) -- cycle;
\draw [rounded corners,line width=1.2,fill opacity=0.25,fill=white!40!black] (-2.4,-0.4) -- (-2.4,2.4) -- (-1.6,2.4) -- (-1.6,-0.4) -- cycle;
\draw [rounded corners,line width=1.2,fill opacity=0.25,fill=white!40!black] (-0.4,2.4) -- (0,2.4) -- (2.4,0) -- (2.4,-0.4) -- (2,-0.4) -- (-0.4,2) -- cycle;
\node at (1.5,-2.05) {\textbf{2}};
\node at (-2,1.5) {\textbf{2}};
\node at (0,2) [below right] {\textbf{2}};

\draw [rounded corners,line width=1.2,fill opacity=0.5,fill=white!40!black] (-1.4,-1.4) -- (-1.4,-0.6) -- (-0.6,-0.6) -- (-0.6,-1.4) -- cycle;
\draw [rounded corners,line width=1.2,fill opacity=0.5,fill=white!40!black] (1.6,-1.4) -- (1.6,-0.6) -- (2.4,-0.6) -- (2.4,-1.4) -- cycle;
\draw [rounded corners,line width=1.2,fill opacity=0.5,fill=white!40!black] (-1.4,1.6) -- (-1.4,2.4) -- (-0.6,2.4) -- (-0.6,1.6) -- cycle;
\node at (-1,-1) [below right] {\textbf{3}};
\node at (2,-1) [below right] {\textbf{3}};
\node at (-1,2) [below right] {\textbf{3}};

\draw [rounded corners,line width=1.2,fill opacity=0.75,fill=white!40!black] (-1.4,-0.4) -- (-1.4,1.3) -- (0,1.3) -- (1.4,0) -- (1.4, -1.4) -- (-0.4,-1.4) -- (-0.4,-0.4) -- cycle;
\node at (-1,1) [below right] {\textbf{4}};

\draw [rounded corners,line width=1.2,fill opacity=0.9,fill=white!40!black] (-0.4,-0.4) -- (-0.4,0.4) -- (0.4,0.4) -- (0.4,-0.4) -- cycle;
\node at (0,0) [below right] {\textbf{5}};

\foreach \x in {-2,-1,...,4}{
  \foreach \y in {-2,-1,...,4}{
    \node[draw,circle,inner sep=0.5pt,fill] at (\x,\y) {};
  }
}

\node at (-2,-2) [below] {\scriptsize$(5,\text{-}1,\text{-}1)$};
\node at (4,-2) [below] {\scriptsize$(\text{-}1,5,\text{-}1)$};
\node at (-2,4) [below] {\scriptsize$(\text{-}1,\text{-}1,5)$};

\end{tikzpicture}}

%% file: toric-fano.bbl
\providecommand{\bysame}{\leavevmode\hbox to3em{\hrulefill}\thinspace}
\providecommand{\MR}{\relax\ifhmode\unskip\space\fi MR }
\providecommand{\MRhref}[2]{%
  \href{http://www.ams.org/mathscinet-getitem?mr=#1}{#2}
}
\providecommand{\href}[2]{#2}
\begin{thebibliography}{BVdV79}

\bibitem[AK77]{altman:77a}
Allen~B. Altman and Steven~L. Kleiman, \emph{Foundations of the theory of
  {F}ano schemes}, Compositio Math. \textbf{34} (1977), no.~1, 3--47.
  \MR{0569043 (58 \#27967)}

\bibitem[BVdV79]{barth:81a}
W.~Barth and A.~Van~de Ven, \emph{Fano varieties of lines on hypersurfaces},
  Arch. Math. (Basel) \textbf{31} (1978/79), no.~1, 96--104. \MR{510081
  (80j:14004)}

\bibitem[CDR08]{casagrande:08a}
Cinzia Casagrande and Sandra Di~Rocco, \emph{Projective {$\Bbb Q$}-factorial
  toric varieties covered by lines}, Commun. Contemp. Math. \textbf{10} (2008),
  no.~3, 363--389. \MR{2417921}

\bibitem[CI15]{ilten:15a}
Melody Chan and Nathan Ilten, \emph{Fano schemes of determinants and
  permanents}, Algebra Number Theory \textbf{9} (2015), no.~3, 629--679.
  \MR{3340547}

\bibitem[CLS11]{CLS}
David~A. Cox, John~B. Little, and Henry~K. Schenck, \emph{Toric varieties},
  Graduate Studies in Mathematics, vol. 124, American Mathematical Society,
  Providence, RI, 2011. \MR{2810322}

\bibitem[DM98]{debarre:98a}
Olivier Debarre and Laurent Manivel, \emph{Sur la vari\'et\'e des espaces
  lin\'eaires contenus dans une intersection compl\`ete}, Math. Ann.
  \textbf{312} (1998), no.~3, 549--574. \MR{1654757 (99j:14048)}

\bibitem[EH00]{eisenbud:00a}
David Eisenbud and Joe Harris, \emph{The geometry of schemes}, Graduate Texts
  in Mathematics, vol. 197, Springer-Verlag, New York, 2000. \MR{1730819}

\bibitem[FI16]{furukawa:16a}
Katsuhisa Furukawa and Atsushi Ito, \emph{A combinatorial description of dual
  defects of toric varieties}, Preprint, 2016.

\bibitem[HMP98]{harris:98a}
Joe Harris, Barry Mazur, and Rahul Pandharipande, \emph{Hypersurfaces of low
  degree}, Duke Math. J. \textbf{95} (1998), no.~1, 125--160. \MR{1646558
  (99j:14043)}

\bibitem[Ilt14]{ilten:14a}
Nathan Ilten, \emph{{F}ano schemes of lines on toric surfaces},
  \href{http://arxiv.org/abs/1411.3025}{\nolinkurl{arXiv:1411.3025}} [math.AG],
  2014.

\bibitem[IT16]{ilten:16a}
Nathan Ilten and Zach Teitler, \emph{Product rank of the $3\times3$ determinant
  and permanent}, Canad. Math. Bull. \textbf{59} (2016), 311--319.

\bibitem[Ito15]{ito:15a}
Atsushi Ito, \emph{Algebro-geometric characterization of {C}ayley polytopes},
  Adv. Math. \textbf{270} (2015), 598--608. \MR{3286545}

\bibitem[Lan97]{langer:97a}
Adrian Langer, \emph{Fano schemes of linear spaces on hypersurfaces},
  Manuscripta Math. \textbf{93} (1997), no.~1, 21--28. \MR{1446187 (98g:14048)}

\bibitem[Lan13]{landsberg}
J.M. Landsberg, \emph{Geometric complexity theory: an introduction for
  geometers}, arXiv:1305.7387v2 [math.AG], 2013.

\bibitem[LK14]{larsen:14a}
Paul Larsen and Franz Kir{\'a}ly, \emph{Fano schemes of generic intersections
  and machine learning}, Internat. J. Algebra Comput. \textbf{24} (2014),
  no.~7, 923--933. \MR{3286146}

\bibitem[LMR13]{landsberg:13a}
Joseph~M. Landsberg, Laurent Manivel, and Nicolas Ressayre, \emph{Hypersurfaces
  with degenerate duals and the geometric complexity theory program}, Comment.
  Math. Helv. \textbf{88} (2013), no.~2, 469--484. \MR{3048194}

\end{thebibliography}
